\documentclass[preprint,3p,times]{elsarticle}

\usepackage{amssymb}
\usepackage{amsmath, amsfonts}
\usepackage{mathtools}
\usepackage{booktabs}
\usepackage{bm}
\usepackage{enumitem}
\usepackage{subcaption}
\usepackage{algorithm}
\usepackage{algpseudocode}

\journal{}

\begin{document}

\begin{frontmatter}



\title{A new operator splitting method for coupled PDE problems, with an application to poromechanics}


\author[label1]{Francesca Marcon} %% Author name
\author[label1]{Stefano Scial\`o}

\affiliation[label1]{organization={Dipartimento di Scienze Matematiche ``G. L. Lagrange'', Politecnico di Torino},%Department and Organization
            addressline={c.so Duca Degli Abruzzi, 24}, 
            city={Turin},
            postcode={10129}, 
            state={},
            country={Italy}}

%% Abstract
\begin{abstract}
The present work proposes a new operator splitting technique. The method leverages ideas used in PDE-constrained domain decomposition, adapting them to the objective of decoupling equations defined on the same geometric domain. The methodology is applied here to the resolution of coupled flow and deformation in the linear quasi-static Biot problem. An overview of its main characteristics is provided, also in comparison with the well established fixed-stress split for poroelasticity.

%The main novelty of the proposed approach lies in the possibility of solving the sub-problems in parallel, in contrast to the other existing strategies that require a sequential resolution. The methodology here presented for poro-elasticity problems, can be extended readily to a broad range of coupled problems. 
\end{abstract}



%% Keywords

\end{frontmatter}


\newcommand{\C}{\bm{C}}
\newcommand{\K}{\frac{\kappa}{\mu_f}}
\newcommand{\bepsilon}{\bm{\epsilon}}
\newcommand{\bu}{\bm{u}}
\newcommand{\bs}[1]{\boldsymbol{#1}}
\newcommand{\spaceVm}{\bs{V}^m}
\newcommand{\spaceVf}{V^f}
\newcommand{\Psiu}{\Psi^{\nabla\cdot \bu}}
\newcommand{\divu}{\nabla\cdot \bu}
\newcommand{\Psip}{\Psi^p}
\newcommand{\lebl}[2][2]{\mathrm{L}^{#1}\ifstrempty{#2}{}{\!(#2)}}

\section{Introduction}
\label{Intro}
The present work proposes a novel operator splitting method, based on PDE constrained optimization. The approach relies on the introduction of additional internal variables to the problem, with the purpose of decoupling the equations in a system. The new variables, indeed, represent a copy of the unknown quantities appearing in more than one equation of the original problem. The problem is thus rewritten as an optimization problem, in which a functional, expressing the mismatch between the new and the original quantities is minimized, subject to the constraints given by the equations of the system. The optimization problem is solved via a gradient based approach in which system equations can be solved independently. The formulation borrows ideas from \cite{BPSa}, where an optimization based domain decomposition approach is developed, and extends its use to a different purpose.

This formulation leads to a splitting of the original problem into smaller sub-problems that can be solved simultaneously, thus enabling the use of parallel computing. Further, independent discretizations can be used for the involved variables. This feature allows to choose the desired approximation level for each of the sub-problems, thus combining efficiency with effectiveness of simulations. 

The proposed approach can be used for general coupled problems. However, for clarity of exposition, it will be presented here for the resolution of the coupled problem of fluid flow and mechanical deformation in poroelastic materials, \cite{Showalter2000}. Numerical resolution strategies for this problem have been proposed by several authors: in \cite{Rodrigo2016} a two-field stabilized formulation can be found; in \cite{Berger2015,Oyarzua2016} two different three-field formulations are presented; a four-field problem is proposed in \cite{Ahmed2019}, whereas a five-field method is analyzed in \cite{Zhao2022}, in which stresses, displacements, fluid pressure and fluid flow appear as unknowns along with an additional variable used to enforce symmetry of the stresses.

Given the large size of the resulting discrete system, techniques exist to decouple the flow subproblem from the mechanics subproblem. The fixed-stress operator splitting \cite{KIM2011CMAME} is a well established alternative to the monolithic solution of the problem. The flow problem is solved first, using the displacement field available from the previous iteration, after which the updated displacement is computed. The convergence of this scheme is studied in \cite{,Mikelic2013,Storvik2019}. The fixed-stress splitting is used as a reference in this work.

For simplicity, here we consider the quasi-static two field formulation of the Biot problem in a domain $\Omega \subset \mathbb{R}^2$ and in a time interval $(0,\, T)$. However, we note that the proposed methodology is independent of the chosen formulation of the problem. We denote by $\partial \Omega$ the boundary of $\Omega$, and further introduce two, possibly different partitions of $\partial \Omega$, defined, respectively, as: $\partial \Omega:=\Gamma^m_D \cup \Gamma^m_N$, with $\Gamma^m_D \cap \Gamma^m_N=\emptyset$ and $|\Gamma^m_D|>0$; and $\partial \Omega:=\Gamma^f_D\cup \Gamma^f_N$, with $\Gamma^f_D\cap \Gamma^f_N=\emptyset$, and $|\Gamma^f_D|>0$.

The following spaces are then defined: $\mathrm{H}^1_{\Gamma^m_D}(\Omega)$ and $\mathrm{H}^1_{\Gamma^f_D}(\Omega)$, as the first order Sobolev spaces of functions having null trace on $\Gamma^m_D$ and $\Gamma^f_D$, respectively. We set, for brevity, $\spaceVm:=[\mathrm{H}^1_{\Gamma_D^m}(\Omega)]^2$, and $\spaceVf:=\mathrm{H}^1_{\Gamma^f_D}(\Omega)$.
Then, the problem of interest reads: \emph{find $\bm{u}\in \lebl{0,T;\spaceVm}$ and $p\in\lebl{0,T;\spaceVf}$:}
\begin{eqnarray}
    &&\int_{\Omega}\C\bepsilon(\bu):\bepsilon(\bm{v}) \, \mathrm{d}\Omega -\alpha \int_{\Omega} p \nabla\cdot \bm{v} \, \mathrm{d}\Omega=\int_\Omega \bm{b}\,{\bm{v}} \, \mathrm{d}\Omega \quad \forall \bm{v}\in \spaceVm \label{mechEq} \\
    && \int \alpha \frac{\partial}{\partial t} \left(\nabla \cdot \bu \right) q \, \mathrm{d}\Omega + \int_\Omega \K \nabla p \nabla q \, \mathrm{d}\Omega = \int_\Omega g q \, \mathrm{d}\Omega \quad \forall q \in \spaceVf \label{flowEq} \\
    && \bs{u}=\bs{u}_0, \ p=p_0,  \qquad \text{in} \ \Omega, \ t=0, \label{initCond}
\end{eqnarray}
where $\bepsilon(\bm{u})=\frac12\left(\nabla \bm{u}+(\nabla \bm{u})^T\right)$, $\C$ is the elasticity tensor, $\bm{b}$ a load, $\alpha>0$ the Biot coefficient, $\kappa$ is the hydraulic permeability of the porous matrix, $\mu_f$ fluid viscosity, and $g$ a source term. We have assumed homogeneous Neumann conditions on $\Gamma^m_N$ and $\Gamma^f_N$, and $\bs{u}_0$, $p_0$ are given functions defined in $\Omega$ at the initial time. We have taken, for simplicity $\frac{1}{M}=0$, with $M$ the Biot modulus.

\section{Parallel operator splitting}
Let us introduce a discretization of the time interval $(0,\, T)$, made of $N^t+1$ time intervals $[t^k,t^{k+1}]$ of equal length $dt$, for $k=0,\ldots,N^t$, and let us denote by $\bu^k$, $p^k$ the unknowns at time $t^k$. By using the backward Euler Method, we obtain the following semi-discrete in time problem: \emph{for $k=1,\ldots,N^t+1$, $dt>0$ and given $\bu^{k-1},\,p^{k-1}$, find $\bm{u}^k\in\spaceVm$ and $p^k\in\spaceVf$:}
\begin{eqnarray}
    &&\int_{\Omega}\C\bepsilon(\bu^k):\bepsilon(\bm{v}) \, \mathrm{d}\Omega -\alpha \int_{\Omega} p^k \nabla\cdot \bm{v} \, \mathrm{d}\Omega=\int_\Omega \bm{b}\,{\bm{v}} \, \mathrm{d}\Omega, \quad \forall \bm{v}\in \spaceVm \label{mechEq_t} \\
    && \frac{\alpha}{dt} \int_\Omega \left(\nabla \cdot \bu^k \right) q \, \mathrm{d}\Omega + \int_\Omega \K \nabla p^k \nabla q \, \mathrm{d}\Omega = \nonumber \\
    && \hspace{3cm}\int_\Omega g q \, \mathrm{d}\Omega +\frac{\alpha}{dt}  \int_\Omega \left(\nabla \cdot \bu^{k-1} \right) q \, \mathrm{d}\Omega, \quad \forall q \in \spaceVf \label{flowEq_t}
\end{eqnarray}

We now introduce the additional unknowns $\Psiu \in \mathrm{L}^2(\Omega)$ and $\Psip \in\mathrm{L}^2(\Omega)$ and rewrite Problem~\eqref{mechEq_t}-\eqref{flowEq_t}, equivalently, as follows: \emph{for $k=1,\ldots,N^t$, $dt>0$ and given $\bu^{k-1},\,p^{k-1}$, find $\bm{u}^k\in \spaceVm$, $p\in\spaceVf$, and $\Psiu \in \mathrm{L}^2(\Omega)$ and $\Psip\in\mathrm{L}^2(\Omega)$:}
\begin{eqnarray}
    &&\min_{\Psiu,\Psip} J:=\frac12 \left(\|(\nabla\cdot \bu^k)-\Psiu\|^2_{\mathrm{L}^2(\Omega)} + \|p^k-\Psip\|^2_{\mathrm{L}^2(\Omega)}\right)\\
    && \text{such that:} \nonumber \\
    && \int_{\Omega}\C\bepsilon(\bu^k):\bepsilon(\bm{v}) \, \mathrm{d}\Omega -\alpha \int_{\Omega} \Psip (\nabla\cdot \bm{v}) \, \mathrm{d}\Omega=\int_\Omega \bm{b}\,{\bm{v}} \, \mathrm{d}\Omega \quad \forall \bm{v}\in \spaceVm \\
    && \frac{\alpha}{dt} \int_\Omega (\Psiu) q  \, \mathrm{d}\Omega + \int_\Omega \K \nabla p \nabla q \, \mathrm{d}\Omega = \nonumber \\
    && \hspace{3cm}\int_\Omega g q  \, \mathrm{d}\Omega + \frac{\alpha}{dt}\int_\Omega \left(\nabla \cdot \bu^{k-1} \right) q \, \mathrm{d}\Omega, \quad \forall q \in \spaceVf
\end{eqnarray}
Clearly $\Psiu$, $\Psip$ depend on $k$, but this dependence is omitted to simplify the notation. We also note that the regularity needed for $\Psip$ is lower than that of its counterpart $p^k$. In addition, we defined $\Psiu$ to match the divergence of the displacements, so a low regularity is also sufficient for this quantity. Different choices are also possible. 

Let us now consider the problem to be solved at each time step, and so, we drop the time index $k$, for simplicity of notation. 
The fully discrete problem follows by introducing space discretizations for all the variables involved. Independent discretizations may be used as well for the involved variables. Thus, we introduce a mesh $\mathcal{T}^m$ with mesh parameter $h_m$ for the displacement field, a mesh $\mathcal{T}^f$ with parameter $h_f$ for the flow field, and meshes $\mathcal{T}^{\divu}$ with parameter $h_{\nabla\cdot \bm{u}}$ and $\mathcal{T}^{p}$ with parameter $h_{p}$ for $\Psiu$ and $\Psip$, respectively. Linear Lagrangian finite elements are used for displacements and for the pressure, while $\Psiu$ and $\Psip$ are approximated by piecewise constant functions, other choices being also possible. We denote by $\spaceVm_h$, $\spaceVf_h$, $Q_h^{\divu}$ and $Q_h^p$ the discrete spaces used for $\bu^k$, $p^k$, $\Psiu$ and $\Psip$, respectively. 
The algebraic system is derived proceeding by collecting into matrices the integrals of the basis functions of the discrete spaces, and into arrays the degrees of freedom of the unknowns.

Let $\{\bs{\varphi}^m_l\}_{l=1}^{\mathcal{N}^m}$ be a set of basis functions of $\spaceVm_h$ with dimension $\mathcal{N}^m$, $\{\varphi_l\}_{l=1}^{\mathcal{N}^f}$ be a set of basis functions of $\spaceVf$ with dimension $\mathcal{N}^f$ and let $\{\vartheta^{\star}_l\}_{l=1}^{\mathcal{N}^{\star}}$ be, for $\star\in\{\divu,p\}$, the sets of basis functions of $Q_h^{\star}$, with dimension $\mathcal{N}^{\star}$. 
The discrete variables are thus written as:
\begin{equation}
    \bs{u}_h = \sum_{l=1}^{\mathcal{N}^m} \mathrm{u}_l \bs{\varphi}^m_l\,,
    \quad 
    p_h = \sum_{l=1}^{\mathcal{N}^f} \limits \mathrm{p}_l^f \varphi_l^f,\;\; 
\quad
\Psi^\star = \sum_{l=1}^{\mathcal{N}^{\star}}\limits
\psi^{\star}_l\vartheta_l^{\star}, \, \star\in\{\divu,p\}
\end{equation}

We then introduce the following matrices: $\bs{K}\in\mathbb{R}^{\mathcal{N}^m\times\mathcal{N}^m}$, $\mathbf{A}\in\mathbb{R}^{\mathcal{N}^f\times \mathcal{N}^f}$, $\mathbf{B}\in\mathbb{R}^{\mathcal{N}^m\times \mathcal{N}^p}$ and $\mathbf{D}\in\mathbb{R}^{\mathcal{N}^f\times \mathcal{N}^{\divu}}$, defined as:
\begin{eqnarray}
    && (\bs{K})_{l,j} = \int_{\Omega} \bs{C} \bs{\epsilon}(\bs{\varphi}^m_j) : \bs{\epsilon}(\bs{\varphi}^m_l), \quad  (\mathbf{A})_{l,j} = \int_\Omega \K \nabla \varphi_l
    \nabla \varphi_j, \\
    && (\bs{B})_{l,j}=\alpha\int_\Omega (\nabla\cdot\bs{\varphi}^m_l)\vartheta^p_j, \quad (\bs{D})_{l,j}=-\frac{\alpha}{dt}\int_\Omega \varphi_l \,\vartheta^{\divu}_j
\end{eqnarray}
and matrices: $\bs{M}^{\divu}\in\mathbb{R}^{\mathcal{N}^m\times\mathcal{N}^m}$, $\bs{M}^{\Psiu}\in\mathbb{R}^{\mathcal{N}^{\divu}\times\mathcal{N}^{\divu}}$, $\bs{E}^{\bu}\in\mathbb{R}^{\mathcal{N}^{\divu}\times\mathcal{N}^{m}}$, $\bs{M}^{p}\in\mathbb{R}^{\mathcal{N}^f\times\mathcal{N}^f}$, $\bs{M}^{\Psip}\in\mathbb{R}^{\mathcal{N}^{p}\times\mathcal{N}^{p}}$, $\bs{E}^{p}\in\mathbb{R}^{\mathcal{N}^p\times\mathcal{N}^{f}}$, defined as:
\begin{eqnarray}
    && (\bs{M}^{\divu})_{l,j} = \int_{\Omega} (\nabla\cdot\bs{\varphi}^m_j)(\nabla\cdot\bs{\varphi}^m_l) , \quad  (\mathbf{M}^{\Psiu})_{l,j} = \int_\Omega \vartheta^{\divu}_j \,\vartheta_l^{\divu}\\
    && (\bs{E}^{\bu})_{l,j}=-\int_\Omega \vartheta^{\divu}_l(\nabla\cdot\bs{\varphi}^m_j), \quad (\bs{M}^p)_{l,j}=\int_\Omega \varphi_j\,\varphi_l \\
    && (\bs{M}^{\Psip})_{l,j} = \int_{\Omega} \vartheta^{p}_j \,\vartheta_l^{p}, \quad  (\mathbf{E}^{p})_{l,j} =- \int_\Omega \vartheta^{p}_l \,\varphi_j.
\end{eqnarray}
We can formulate the algebraic version of the problem that is solved at each time step:
\begin{eqnarray}
    && \min_{\bs{\psi}^{\divu},\bs{\psi}^{p}}\frac{\eta}{2}\left( \bs{\mathrm{u}}^T\bs{M}^{\divu} \bs{\mathrm{u}} + (\bs{\psi}^{\divu})^T \bs{E}^{\divu}\bs{\mathrm{u}} + \bs{\mathrm{u}}^T (\bs{E}^{\divu})^T \bs{\psi}^{\divu}+ (\bs{\psi}^{\divu})^T\bs{M}^{\Psip}\bs{\psi}^{\divu}\right)+\nonumber \\
    && \hspace{1cm}\frac12\left(\bs{\mathrm{p}}^T\bs{M}^{p} \bs{\mathrm{p}} + (\bs{\psi}^{p})^T \bs{E}^{p}\bs{\mathrm{p}} + \bs{\mathrm{p}}^T (\bs{E}^{p})^T \bs{\psi}^{p}+(\bs{\psi}^p)^T\bs{M}^{\Psip}\bs{\psi}^p\right)\\
    && \text{such that} \nonumber \\
    && \bs{K} \bs{\mathrm{u}} - \bs{B}\bs{\psi}^p = \bs{\mathrm{b}} \label{mechProb_discr}\\
    && \bs{A} \bs{\mathrm{p}} - \bs{D}\bs{\psi}^{\divu} =\bs{\mathrm{g}} \label{flowProb_discr}
\end{eqnarray}
where $\bs{\mathrm{u}}$, $\bs{\mathrm{p}}$, $\bs{\psi}^{\divu}$, and $\bs{\psi}^p$ are column arrays of degrees of freedom, while $\bs{\mathrm{b}}$ and $\bs{\mathrm{g}}$ collect all known terms.
The coefficient $\eta$ is a scaling term, that, in the discrete problem can be used to balance the two terms in the functional, helping its minimization. 
%It can be chosen at each time iteration $k$ as $\eta=\frac{\|\bs{\mathrm{p}}^{k-1}\|^2}{\|\bs{\mathrm{u}}^{k-1}\|^2}$.

The above constrained minimization problem can be recast in an unconstrained minimization problem, by exploiting the linearity of the constraints, yielding:
\begin{eqnarray}
    && \min_{{\bs{\psi}^{\divu},\bs{\psi}^{p}}} \frac12\left((\bs{\psi}^{\divu})^T \bs{\mathcal{G}^1} \bs{\psi}^{\divu}+(\bs{\psi}^{p})^T \bs{\mathcal{G}^2} \bs{\psi}^{p}+ (\bs{\psi}^{\divu})^T \bs{\mathcal{H}}\bs{\psi}^p + \right. \nonumber \\
    && \left.\hspace{3cm} (\bs{\psi}^{p})^T \bs{\mathcal{H}}^T\bs{\psi}^{\divu}+
    (\bs{\psi}^{\divu})^T\bs{\tilde{\mathrm{b}}}+(\bs{\psi}^{p})^T\bs{\tilde{\mathrm{g}}}\right)
    \label{eq:unconstrOpt}
\end{eqnarray}
with:
\begin{displaymath}
    \bs{\mathcal{G}^1}=\eta\bs{B}^T\bs{K}^{-T}\bs{M}^{\divu}\bs{K}^{-1}\bs{B}+\eta\bs{M}^{\Psiu},  \quad
    \bs{\mathcal{G}^2}=\bs{D}^T\bs{A}^{-T}\bs{M}^{p}\bs{A}^{-1}\bs{D}+\bs{M}^{\Psip},
\end{displaymath}
\begin{displaymath}
    \bs{\mathcal{H}}=\eta\bs{E}^{\divu}\bs{K}^{-1}\bs{B}+\bs{D}^T \bs{A}^{-T}(\bs{E}^p)^T,  \quad \bs{\tilde{\mathrm{b}}}=\eta\bs{E}^{\divu}\bs{K}^{-1}\bs{\mathrm{b}}, \quad \bs{\tilde{\mathrm{g}}}=\bs{E}^{p}\bs{A}^{-1}\bs{\mathrm{g}},
\end{displaymath}
and constant terms are omitted. 
This problem is solved via a gradient based approach, and the gradient direction at point $(\bs{\psi}^{\divu})_0$, $(\bs{\psi}^{p})_0$  is computed via Algorithm~\ref{alg:initial-direction}.

% \begin{enumerate}[itemsep=0.8em]
%     \item compute $(\bs{\mathrm{u}})_0$, $(\bs{\mathrm{p}})_0$ solving \eqref{mechProb_discr}-\eqref{flowProb_discr};
%     \item compute $(\bs{\lambda^u})_0$, $(\bs{\lambda^p})_0$ solving the dual problems  
%     $$\bs{K}^T(\bs{\lambda^u})_0=\eta\left(\bs{M}^{\divu}(\bs{\mathrm{u}})_0+(\bs{E}^{\divu})^T(\bs{\psi}^{\divu})_0\right),$$ $$\bs{A}^T(\bs{\lambda^p})_0=\bs{M}^{p}(\bs{\mathrm{p}})_0+(\bs{E}^p)^T(\bs{\psi}^{p})_0;$$
%     \item set
%     \begin{equation}
%         \bs{d}_0=\begin{bmatrix}
%             \bs{D}^T(\bs{\lambda^p})_0+\eta\bs{E}^{\divu}(\bs{\mathrm{u}})_0+\eta\bs{M}^{\Psiu}(\bs{\psi}^{\divu})_0 \\
%             \bs{B}^T(\bs{\lambda^u})_0+\bs{E}^{p}(\bs{\mathrm{p}})_0+\bs{M}^{\Psip}(\bs{\psi}^{p})_0
%         \end{bmatrix}.
%         \label{eq:gradDir}
%     \end{equation}
% \end{enumerate}
\begin{algorithm}[ht]
    \caption{Computation of the initial search direction}
    \label{alg:initial-direction}

    \begin{algorithmic}[1]

        \State Compute $(\bs{\mathrm{u}})_0$ and
        $(\bs{\mathrm{p}})_0$ by solving
        \eqref{mechProb_discr}--\eqref{flowProb_discr}. \label{Algo:1}

        \State Compute $(\bs{\lambda^u})_0$ and
        $(\bs{\lambda^p})_0$ by solving the dual problems \label{Algo:2}
        \begin{eqnarray}
                \bs{K}^{T}(\bs{\lambda^u})_0
                &=&
                \eta\left(
                    \bs{M}^{\divu}(\bs{\mathrm{u}})_0
                    +(\bs{E}^{\divu})^{T}
                    (\bs{\psi}^{\divu})_0
                \right),                                              \label{eq:dualU}  \\
                \bs{A}^{T}(\bs{\lambda^p})_0
                &=&
                \bs{M}^{p}(\bs{\mathrm{p}})_0
                +(\bs{E}^{p})^{T}(\bs{\psi}^{p})_0. \label{eq:dualP}
        \end{eqnarray}

        \State Set the initial search direction as
        \begin{equation} 
            \bs{d}_0 =
            \begin{bmatrix}
                \bs{D}^{T}(\bs{\lambda^p})_0
                +\eta\bs{E}^{\divu}(\bs{\mathrm{u}})_0
                +\eta\bs{M}^{\Psiu}(\bs{\psi}^{\divu})_0
                \\[2mm]
                \bs{B}^{T}(\bs{\lambda^u})_0
                +\bs{E}^{p}(\bs{\mathrm{p}})_0
                +\bs{M}^{\Psip}(\bs{\psi}^{p})_0
            \end{bmatrix}.
            \label{eq:gradDir}
        \end{equation} \label{Algo:3}

    \end{algorithmic}
\end{algorithm}
It is noted that the descent direction can be computed by solving the primal and dual problems independently, and this can be performed in parallel. For self-adjoint operators, as the ones considered in this work, primal and dual problems coincide. Moreover, often, very efficient solver exist for the decoupled problems, which makes the splitting even more convenient.

\subsection*{Fixed stress split}
For reference we report here the formulation of the fixed-stress scheme used for comparison in the numerical result Section.

Let $k>0$ be the time-step index, $i\geq 1$ be the iteration index and assume $\bs{u}_h^{k-1}$, $(p_h)^{k-1} $, and $\bs{u}_h^{k,i-1}$, $(p_h)^{k,i-1} $ are given.

First find  $(p_h)^{k,i}\in V_h^{f}$ given by:
    \begin{eqnarray*}
&&\frac{L}{dt}\int_{\Omega}  p_h^{k,i}  q_h 
+ \int_{\Omega} \K \nabla p_h^{k,i} \, \nabla q_h = \frac{\alpha}{d t} \int_{\Omega} \nabla \cdot \bs{u}_h^{k-1}  q_h+ \int_{\Omega} g\, q_h  
\\
&& \hspace{1cm}- \frac{\alpha}{dt}\int_{\Omega}  \nabla\cdot \bs{u}_h^{k,i-1}  q_h 
+\frac{L}{dt}\int_{\Omega}p_h^{k,i-1}  q_h 
\quad \forall  q_h\in V_h^{f} 
\end{eqnarray*}
Then, find $\bs{u}_h^{k,i}\in\bs{V}^m_h$ such that:
    \begin{eqnarray*}
 \int_{\Omega} \bs{C} \bs{\epsilon}(\bs{u}^{k,i}_h) : \bs{\epsilon}(\bs{v}_h)= \int_{\Omega} \alpha p_h^{k,i} \nabla\cdot \bs{v}_h
\quad \forall \; \bs{v}_h\in \bs{V}_h^{m},
\end{eqnarray*}
and $L=\frac{\alpha^2}{E}$, as in \cite{KIM2011CMAME}, being $E$ the drained bulk modulus.

\section{Numerical Results}

As an application, we consider the Terzaghi problem \cite{terzaghi,Kadeethum2020} in a homogeneous square domain $\Omega$ with edge length of $1 \, \mathrm{m}$. For the displacement problem, a compressive traction of $-1 \,\mathrm{kN/m}$ is set on the top edge $y=1$, while only vertical displacement is allowed on the lateral edges, and clamping is imposed on the bottom edge. For the flow problem, a homogeneous Dirichlet condition is prescribed on the edge $y=1$, and homogeneous Neumann conditions on the rest of the boundary. %Domain and boundary conditions are sketched in Figure~\ref{fig:domain}.

% \begin{figure}
%     \centering
%     \begin{minipage}{0.47\textwidth}
%     \includegraphics[width=0.80\linewidth]{Figure/Terzaghi.pdf}
%     \caption{Domain and boundary conditions}
%     \label{fig:domain}    
%     \end{minipage}
%     \begin{minipage}{0.04\textwidth}
%         ~
%     \end{minipage}
%     \begin{minipage}{0.47\textwidth}
%     \includegraphics[width=0.99\linewidth]{Figure/figure_11.pdf}
%     \caption{Solution obtained with the proposed approach at different times}
%     \label{fig:ExampleSol}    
%     \end{minipage}
% \end{figure}

For this problem the exact pressure solution is available, given at each time $t>0$ for $x=0.5$ and $0\leq y\leq 1$, by:
\begin{equation}
    p(y,t) = \sum_{m=0}^{\infty} \frac{4}{\pi(2m+1)} \sin\left(\frac{(2m+1)\pi}{2}y\right) \exp{\left(-\frac{(2m+1)^2\pi^2}{4}C_v t\right) }\,,
    \label{eq:pex}
\end{equation}
being
$$C_v = 3E\,\frac{1-\nu}{1+\nu} \K,$$
$E$ the bulk modulus, and $\nu$ the Poisson ratio of the matrix.
%We recall the following relations to compute the Lamé coefficients given the bulk modulus  and the Poisson ratio
%\begin{equation}
%    \lambda = \frac{3E\nu}{1+\nu}\,, \quad \mu = \frac{3E(1-2\nu)}{2(1+\nu)}\,.
%\end{equation}
The following values are assigned: $E=10^3\,\mathrm{kPa}$, $\nu=0.25$, $\kappa = 10^{-12}\,\mathrm{m}^2$, $\mu_f=10^{-6} \,\mathrm{kPa\, s}$, $\alpha = 1$ and $g=0$. 

Four, possibly independent, triangular meshes are defined in $\Omega$ to numerically resolve the problem with the proposed approach. Quantities $h_m$, $h_f$, $h_{\divu}$ and $h_p$ denote the maximum area of the elements in each mesh, as described in the previous Section, and they are collected in the array $h=[h_m,\,h_f,\, h_{\divu},\, h_p]$. The syntax $h=\text{value}$ will be used to denote that the same mesh parameter is used for all discretizations. The time step for the backward Euler method is $dt=1 \, \mathrm{s}$. At time $t=0$ it is $p=0$ and $\bs{u}=0$ in the whole $\Omega$.
The preconditioned conjugate-gradient method is used to solve the optimization problem~\eqref{eq:unconstrOpt}, with scaling parameter $\eta=10^{8}$. Gradient direction $\bs{d}$ of equation~\eqref{eq:gradDir} is corrected using the block-diagonal preconditioner $\bs{\mathcal{P}}$ obtained taking the main diagonals of matrices $\eta \bs{M}^{\Psiu}$ and $\bs{M}^{\Psip}$. 

Let us define the following error:
\begin{equation}
    \text{err}_p^k=\max_{j=1,\ldots,20}|p^{6}(y_j,t^k)-p_h^k(y_j)|,
    \label{eq:discrErr}
\end{equation}
being $0=y_1<y_2<\cdots<y_{19}<y_{20}=1$ equally spaced nodes on the $x=0.5$ line, and $p^6(y,t)$ the function obtained considering the first six terms in the sum of equation~\eqref{eq:pex}.

\begin{figure}
    \begin{subfigure}{0.47\textwidth}
    \includegraphics[width=0.99\linewidth]{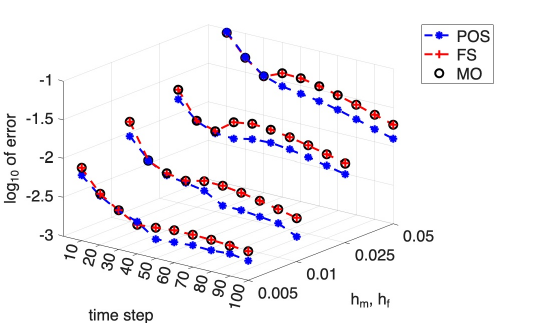}
    \caption{Error $\text{err}_p^k$ (logarithm) for POS and FS at selected time-steps, with error of monolithic approach (MO)}
    \label{subfig0:ErrorComparison}    
    \end{subfigure}
    \begin{subfigure}{0.47\textwidth}
    \includegraphics[width=0.99\linewidth]{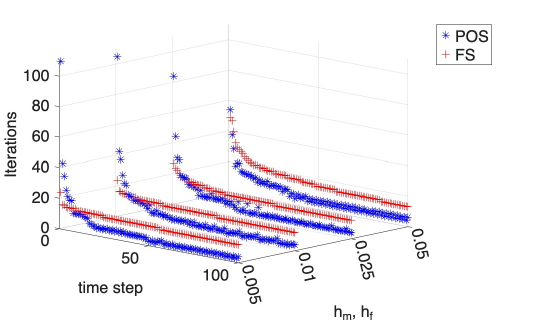}
    \caption{Number of PCG and fixed-stress iterations at each time-step}
    \label{subfig0:IterComparison}    
    \end{subfigure}
    \caption{Comparison of POS and FS for various meshes: $h\in\{0.05, \, 0.025, \, 0.01, \, 0.005\}$}
\end{figure}

\begin{table}
    \centering
    \caption{Computational performance of POS and FS for different values of $h$.}
    \label{tab:performance}
    
    \begin{tabular}{@{}lcccc@{}}
        \toprule
        & \multicolumn{4}{c}{Mesh size} \\
        \cmidrule(lr){2-5}
        Method
        & $h=0.05$
        & $h=0.025$
        & $h=0.01$
        & $h=0.005$ \\
        \midrule
        
        \multicolumn{5}{l}{
            \textit{Average time per iteration} 
            ($\times 10^{-4}$ $\mathrm{[s]}$)
        } \\
        \addlinespace[2pt]
        POS & 0.35 & 0.40 & 0.98 & 2.51 \\
        FS  & 0.18 & 0.31 & 0.49 & 1.20 \\
        
        \midrule

        \multicolumn{5}{l}{
            \textit{Pre-processing time} 
            ($\mathrm{[s]}$)
        } \\
        \addlinespace[2pt]
        POS & 0.024 & 0.047 & 0.110 & 0.185 \\
        FS  & 0.026 & 0.030 & 0.043 & 0.064 \\
        
        \midrule
        
        \multicolumn{5}{l}{
            \textit{Total number of iterations}
        } \\
        \addlinespace[2pt]
        POS & 824  & 834  & 798  & 729  \\
        FS  & 1750 & 1367 & 1308 & 1299 \\
        
        \bottomrule
    \end{tabular}
\end{table}

We propose a comparison between the proposed approach, labeled as POS, and the fixed-stress split, labeled as FS. Given the different nature of the two algorithms, a different stopping criterion is selected. For the PCG iterations of the POS method, the criterion is based, as usual, on the norm of the residual relative to the norm of the initial residual. For the FS method, the stopping criterion is based instead on the relative distance of the solution across two iterations in euclidean norm: we set $\text{dist}_{\bu}=\|\bs{\mathrm{u}}^{k,i}-\bs{\mathrm{u}}^{k,i-1}\|_2/\|\bs{\mathrm{u}}^{k,i}\|_2$ and $\text{dist}_{p}=\|\mathrm{p}^{k,i}-\mathrm{p}^{k,i-1}\|_2/\|\mathrm{p}^{k,i}\|_2$, where, we recall that $\bs{\mathrm{u}}^{k,i}$ and $\mathrm{p}^{k,i}$ are the array of unknowns for displacement and pressure, respectively, at time $t^k$ and at the $i$-th fixed-stress iteration. Fixed stress iterations are stopped when both distances fall below a prescribed tolerance. The threshold FS is set at $10^{-6}$, corresponding to the larger value for which the error in \eqref{eq:discrErr} evaluated for the fixed-stress pressure solution overlaps with  the error given by the solution obtained with a monolithic approach, for all the meshes considered. We recall that the so-called monolithic approach consists in solving the coupled poroelasticity system. The tolerance for POS is set to $10^{-4}$. In  Figure~\ref{subfig0:ErrorComparison} the logarithm of error $\text{err}_p^k$ is reported for the pressure solution obtained with the POS, FS and monolithic (MO) approaches at $t^k$, for $k=10,20,\ldots,100$, on four meshes, with parameter $h\in\{0.05, \, 0.025,\,0.01,\,0.005\}$. The same mesh is used for all the methods and, in the case of POS for all the spatial triangulations. It can be seen that the curves for FS and MO are overlapped, and the curve for POS lies slightly below. Figure~\ref{subfig0:IterComparison} shows the number of iteration, for each time-step and for the same meshes required by the POS and FS methods to reach convergence. It can be seen that for all the considered meshes, the POS approach requires more iterations than FS only during the first few time-steps; thereafter it converges in only a few iterations per step. In contrast, the number of FS iterations remains quite constant throughout time advancing. This gives, overall, the total iteration count reported in Table~\ref{tab:performance}. In the same table the average time per iteration and the preprocessing time are also reported for the two approaches on the same considered meshes. As expected, the average time per iteration of the POS approach is about twice that of the FS method, since the computation of the gradient via Algorithm~\ref{alg:initial-direction} requires the solution of the primal and dual linear systems (steps \ref{Algo:1}-\ref{Algo:2}). Also, the pre-processing time is larger for the POS approach, for the computation of the meshes and of the matrices related to the additional variables introduced by the method. However, in consideration of the lower iteration count, the POS method can become computationally convenient for long integration times. Further, the resolution of the primal and dual linear systems \eqref{mechProb_discr}-\eqref{flowProb_discr} and \eqref{eq:dualU}-\eqref{eq:dualP} can be performed in parallel, as well as the pre-processing part. 

\begin{figure}
    \begin{subfigure}{0.47\textwidth}
    \includegraphics[width=0.95\linewidth]{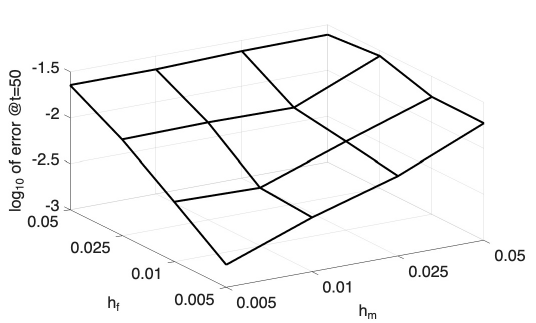}
    \caption{Error (logarithm) at $t=50$}
    \label{subfig1:Err}    
    \end{subfigure}
    \begin{subfigure}{0.47\textwidth}
    \includegraphics[width=0.95\linewidth]{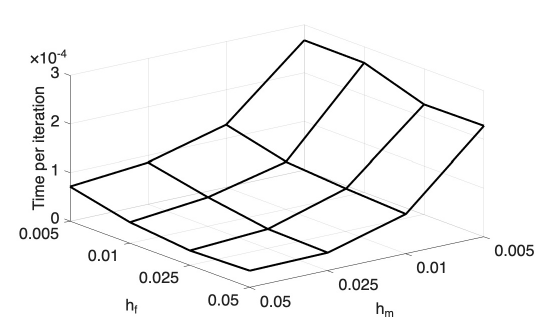}
    \caption{Average time per PCG iteration}
    \label{subfig1:PCGtime}    
    \end{subfigure}
    \begin{subfigure}{0.47\textwidth}
    \includegraphics[width=0.95\linewidth]{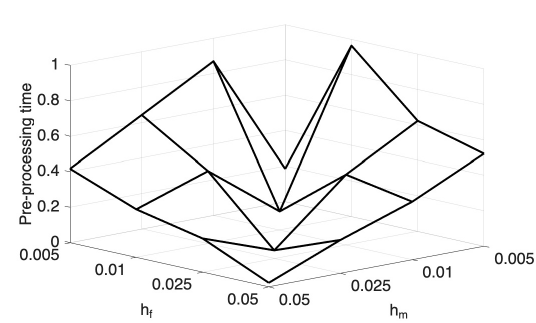}
    \caption{Pre-processing time}
    \label{subfig1:TotalTime}    
    \end{subfigure}
    \begin{subfigure}{0.47\textwidth}
    \includegraphics[width=0.95\linewidth]{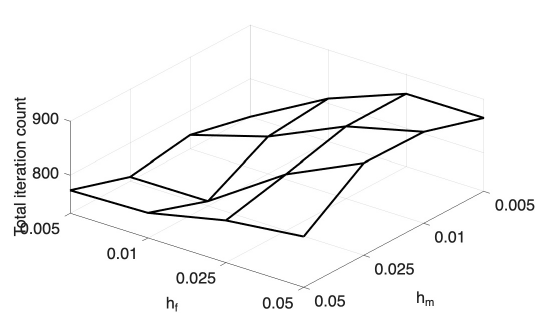}
    \caption{Total iteration count}
    \label{subfig1:TotalIter}    
    \end{subfigure}
    \caption{Behavior of the method at varying of mesh sizes $h_{m}$ and $h_f$}
    \label{fig:BehaviorVsMesh}
\end{figure}
\begin{figure}
    \centering
    \begin{subfigure}{0.47\textwidth}
    \includegraphics[width=0.90\linewidth]{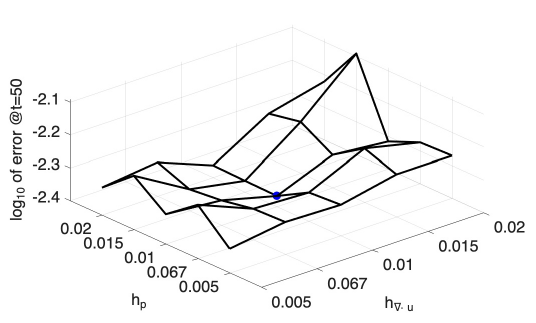}
    \caption{Error (logarithm) at $t=50$}
    \label{subfig2:Err}    
    \end{subfigure}
    \begin{subfigure}{0.47\textwidth}
    \includegraphics[width=0.99\linewidth]{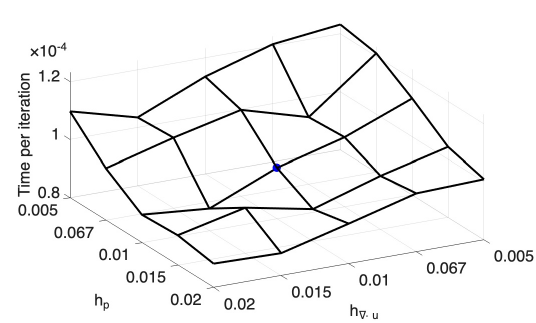}
    \caption{Average time per PCG iteration}
    \label{subfig2:PCGtime}    
    \end{subfigure}
    \begin{subfigure}{0.47\textwidth}
    \includegraphics[width=0.90\linewidth]{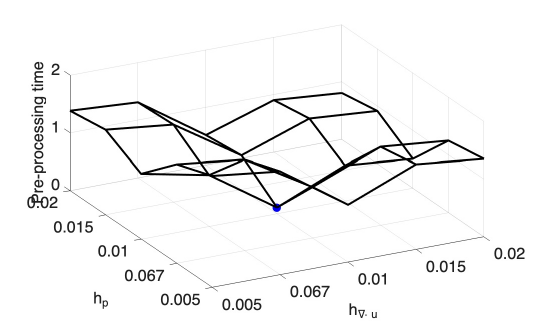}
    \caption{Pre-processing time}
    \label{subfig2:totTime}    
    \end{subfigure}
    \begin{subfigure}{0.47\textwidth}
    \includegraphics[width=0.99\linewidth]{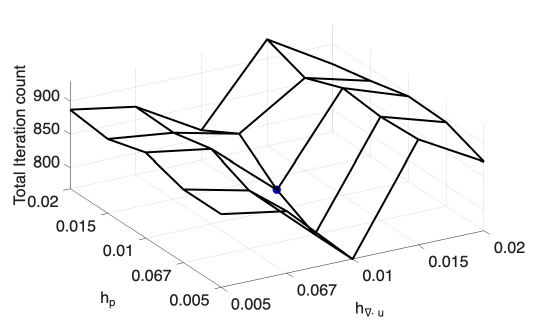}
    \caption{Total iteration count}
    \label{subfig2:totIter}    
    \end{subfigure}
    \caption{Behavior of the method at varying of mesh sizes $h_{\divu}$ and $h_p$, for $h_m=h_f=0.01$.}
    \label{fig:BehaviorVsIntMesh}
\end{figure}

Figure~\ref{fig:BehaviorVsMesh} reports information on the behavior of the method when meshes $\mathcal{T}^m$ and $\mathcal{T}^f$ are different. In all these graphs also the meshes for the auxiliary variables vary, but we always set $\mathcal{T}^{\divu}$ equal to $\mathcal{T}^m$ and $\mathcal{T}^p$ equal to $\mathcal{T}^f$. Panel~\ref{subfig1:Err} shows the logarithm of quantity $\text{err}_p^{50}$, Panel~\ref{subfig1:PCGtime} the average time per iteration, Panel~\ref{subfig1:TotalTime} pre processing time and Panel~\ref{subfig1:TotalIter} the total iteration count for the various mesh combinations. Inspection of panels~\ref{subfig1:Err}-\ref{subfig1:PCGtime} shows that the POS method enables intermediate levels of accuracy and computational cost, by selecting different mesh combinations. Panel~\ref{subfig1:TotalTime} clearly shows that the pre-processing time is significantly lower when $h_m=h_f=h_{\divu}=h_p$. This is directly related to the fact that, in this circumstance, it is not necessary to compute integrals of basis functions defined on different meshes. However the increase in pre-processing time remains limited. Finally, Panel~\ref{subfig1:TotalIter} reveals that the total iteration count does not significantly vary for different mesh combinations.

Figure~\ref{fig:BehaviorVsIntMesh} presents the same quantities for varying $h_{\divu}$ and $h_p$, with $h_m=h_f=0.01$. The largest variations are observed in Panel~\ref{subfig2:totTime}, in the preprocessing time, which increases significantly when $h_{\divu}$, $h_p$ differ from $h_m=h_f$. Time per iteration, shown in Panel~\ref{subfig2:PCGtime} increases marginally with refinement of $h_{\divu}$ and $h_p$, as expected, as they do not influence the dimensions of the linear systems. Moreover error and total iteration count, reported in Panels~\ref{subfig2:Err} and \ref{subfig2:totIter} are only slightly affected by variations of these parameters. It is however important to note that these parameters are expected to have a larger impact in the parallel implementation of the algorithm, as they are directly related to the amount of data to be shared by parallel processes. 

\section{Conclusions}

The present work has introduced a new operator splitting method, suitable for a large variety of applications. The method is based on a reformulation of the original problem as an optimization problem. The method is applied here to the solution of the Biot problem for linear, quasi-static, poroelasticity and is compared with the well-established fixed-stress split for this class of problems. This preliminary investigation shows the applicability of the method, and provides an overview of its characteristics. In particular, the method shows a faster convergence in terms of number of iterations for long-time simulations. Also the method allows different meshes to be used for different unknowns, thus giving a great flexibility in approximation choices. The method is also suitable for parallel implementation.  

\bibliographystyle{elsarticle-num}
\bibliography{bibliography}
\end{document}